\documentclass[12 pt]{amsart2000}
\usepackage{amsmath2000,amsthm2000,amsfonts,amssymb,verbatim,graphics}
\usepackage[all]{xy}
\def \( {\left( }
\def\){\right)}

\def\cb{{\mathcal B}}
\def\cc{{\mathcal C}}

\def\ce{{\mathcal E}}
\def\cf{{\mathcal F}}

\def\cm{{\mathcal M}}

\def\co{{\mathcal O}}
\def\cp{{\mathcal P}}

\def\RR{\mathbb R}

\def\ri{\hat\mu}

\theoremstyle{plain}

\newtheorem{thm}{Theorem}
\newtheorem*{thm*}{Theorem}
\newtheorem{lem}{Lemma}
\newtheorem{cor}{Corollary}
\newtheorem*{cor*}{Corollary}

\newtheorem*{prop*}{Proposition}

\theoremstyle{definition}

\newtheorem{defn}[thm]{Definition}
\newtheorem*{defn*}{Definition}

\newtheorem{exmp*}{Example}

\newtheorem{exmps*}{Examples}

\theoremstyle{remark}

\newtheorem*{rem*}{Remark}

\newtheorem*{rems*}{Remarks}

\newtheorem*{ack*}{Acknowledgment}

\newtheorem{ex}{Example} 

\newtheorem{question}{Question}

\newcommand{\E}{\mathbb{E}}
\renewcommand{\P}{\mathcal{P}}
\newcommand{\B}{\mathcal{B}}
\newcommand{\one}{\mathbf{1}}
\newcommand{\F}{\mathcal{F}}
\newcommand{\Q}{\mathcal{Q}}
\newcommand{\Z}{\mathbb{Z}}

\newcommand{\PP}{{\mathcal{P}}}

\begin{document}

\title
{Measures of maximal relative entropy}

\author{Karl Petersen}
\address{Department of Mathematics,
CB 3250, Phillips Hall,
         University of North Carolina,
Chapel Hill, NC 27599 USA}
\email{petersen@math.unc.edu}

\author{Anthony Quas}
\address{Department of Mathematical Sciences,
University of Memphis,
Memphis, TN 38152-6429}
\email{quasa@msci.memphis.edu}

\author{Sujin Shin}
\address{Department of Mathematics and Statistics,
University of Victoria, 
Victoria, BC V8W 3P4, Canada. 
Department of Mathematics, Ajou University, 
Suwon 422-749,  South Korea}
\curraddr{Department of Mathematics, Korea Advanced Institute of Science
and Technology, Daejon, 305-701, South Korea}
\email{sjs@math.kaist.ac.kr}
\begin{abstract}
Given an irreducible subshift of finite type $X$, a subshift $Y$, a 
factor map $\pi : X \to Y$, 
and an ergodic invariant measure $\nu$ on $Y$, there can exist more
than one ergodic measure on $X$ which projects to $\nu$ and has
maximal entropy among all measures in the fiber, but there is an
explicit bound
on the number of such maximal entropy preimages.
\end{abstract}
\maketitle
\section{Introduction}

It is a well-known result of Shannon and Parry \cite{Shannon, Parry1}
that every irreducible subshift of finite type (SFT) $X$ on a finite alphabet has
a unique measure $\mu_X$ of maximal entropy for the shift
transformation $\sigma$. The maximal measure is
Markov, and its initial distribution and transition probabilities are
given explicitly in terms of the maximum eigenvalue and corresponding
eigenvectors of the 0,1 transition matrix for the subshift. We are
interested in any possible relative version of this result: given
an irreducible SFT $X$, a subshift $Y$, a factor map $\pi : X \to Y$,
and an ergodic invariant measure $\nu$ on $Y$, how many ergodic
invariant measures can there be on $X$ that project under $\pi$ to
$\nu$ and have maximal entropy in the fiber $\pi^{-1}\{ \nu\}$?  We
will show that there can be more than one such ergodic relatively
maximal measure over a given $\nu$, but there are only finitely many. 
In fact, if $\pi$ is a 1-block map, there can be no more than the cardinality
of the alphabet of $X$ (see Corollary \ref{bound}, below). 
Call a measure $\nu$ on $Y$ {\em
  $\pi$-determinate} in case it has a unique preimage of maximal
entropy. We provide some sufficient conditions for $\pi$-determinacy
and give examples of situations in which relatively maximal measures
can be constructed explicitly.

Throughout the paper, unless stated otherwise $X$ will denote an
irreducible SFT, $Y$ a subshift on a finite alphabet, and $\pi : X \to
Y$ a factor map (one-to-one, onto, shift-commuting map). 
By recoding if necessary, we may assume that $X$ is a 1-step
SFT, so that it consists of all (2-sided) sequences on a finite
alphabet consistent with the allowed transitions described by a
directed graph with vertex set equal to the alphabet, and 
that $\pi$ is a 1-block map.
In the following, ``measure'' means ``Borel probability measure'',  $\cc(X)$ denotes the set of
continuous real-valued functions on $X$, $\cm(X)$ the space of
$\sigma$-invariant measures on $X$, and $\ce(X)
\subset \cm(X)$ the set of ergodic measures on $X$.

Some of the interest of this problem 
arises from its connections (discussed in \cite{Petersen00}) with
information-compressing channels \cite{Marcus-PW}, non-Markov
functions of Markov chains \cite{Blackwell, Boyle-T, Burke-R,
  Marcus-PW}, measures of maximal Hausdorff dimension
and measures that maximize, for a given $\alpha >0$, 
the weighted entropy functional 
\begin{equation}
\phi_\alpha (\mu)=\frac{1}{\alpha+1}[h(\mu) + \alpha h(\pi \mu)]
\end{equation}
\cite{Gatzouras-P, Shin, Shin2}, and relative pressure and relative
equilibrium states \cite{Ledrappier-Y, Walters}. 
The theory of pressure and equilibrium states (see \cite{Ruelle,
  Israel, Keller}), relative pressure and relative equilibrium states
\cite{Ledrappier-W, Walters}, and compensation functions
\cite{Boyle-T, Walters} provides basic tools in this area. 
For a factor map $\pi : X \to Y$ between compact topological dynamical systems
and potential function $V \in \cc (X)$, 
Ledrappier and Walters \cite{Ledrappier-W} defined the {\em relative
  pressure} 
$P(\pi, V) : Y \to \RR$ (a Borel measurable function)
and proved
a relative variational principle: 
For each $\nu \in \cm (Y)$,
\begin{equation}
\int_Y P(\pi ,V) \, d \nu = \sup \{h_\mu (X|Y) + \int_X V \, d \mu : \mu
\in \pi^{-1}\nu \} .
\end{equation}
Any measure $\mu$ that attains the supremum is called a {\em relative
  equilibrium state}. 
A consequence is that
the ergodic
measures $\mu$ that have maximal entropy among all measures in 
$\pi^{-1}\{\nu\}$ have
relative entropy given by
\begin{equation}
h_\mu(X|Y)= \int_Y \lim_{n \to \infty}\frac{1}{n}\log
|\pi^{-1}[y_0\dots y_{n-1}]|\, d\nu (y) .
\end{equation}
($|\pi^{-1}[y_0\dots y_{n-1}]|$ is the number of $n$-blocks in $X$
that map under $\pi$ to the $n$-block $y_0\dots y_{n-1}$.)
By the Subadditive Ergodic Theorem, the limit inside the integral
exists a.e.
with respect 
to each ergodic measure $\nu$ on $Y$, and  
it is constant a.e..
The quantity
\begin{equation}
P(\pi ,0)(y)=\limsup_{n \to \infty}\frac{1}{n}\log
|\pi^{-1}[y_0\dots y_{n-1}]|
\end{equation}
is the {\em relative pressure} of the function $0$ over $y \in Y$. 
The maximum possible relative entropy may be thought of as a
``relative topological entropy over $\nu$''; we denote it by
$h_{\text{top}}(X|\nu)$. 

To understand when a Markov measure on $Y$ has a Markov measure on $X$ 
in its preimage under $\pi$, Boyle and Tuncel introduced the idea of a 
compensation function \cite{Boyle-T}, and the concept was developed further 
by Walters \cite{Walters}. Given a factor map $\pi:X \to Y$ between
topological dynamical systems, a {\em
  compensation function} is a continuous function $F:X \to \mathbb R$
such that
\begin{equation}
P_Y(V) = P_X(V \circ \pi + F) \quad\text{for all } V \in \cc (Y) .
\end{equation}
The idea is that, because $\pi : \cm (X) \to \cm (Y)$ is many-to-one,
we always have 
\begin{align}
P_Y(V) &= \sup\{ h_\nu(\sigma) + \int_YV\,d\nu : \nu \in \cm (Y)
  \} \\
 &\leq  \sup\{ h_\mu(\sigma) + \int_X V \circ \pi \,d\mu : \mu \in \cm (X)
   \} ,
\end{align}
 and a compensation function $F$ can take into account, 
for all potential functions $V$ on $Y$ at once, the extra freedom,
information, or free energy that is available in $X$ as compared to
$Y$ because of the ability to move around in fibers over points of
$Y$. 
A compensation function of the form $G \circ \pi$ with $G \in \cc (Y)$
is said to be {\em saturated}.

The machinery of relative equilibrium states and compensation functions
is used to establish the following basic result about relatively maximal measures
\cite{Shin, Walters}:

{\em Suppose that $\nu \in \ce (Y)$ and $\pi \mu = \nu$. Then $\mu$ is
relatively maximal over $\nu$ if and only if there is $V \in \cc (Y)$
such that $\mu$ is an equilibrium state of $V \circ \pi$.}

Notice that if there is a {\em locally constant} saturated compensation
function $G \circ \pi$, then every Markov measure on $Y$ is
$\pi$-determinate 
with Markov relatively maximal lift, because in \cite{Walters} it is
shown that if there is a saturated compensation function $G \circ
\pi$, then 
the relatively maximal measures over an
equilibrium state of $V \in \cc (Y)$ are the equilibrium states of $V
\circ \pi + G \circ \pi$. 

Further, $\mu_X$ is the unique
equilibrium state of the potential function 0 on $X$, the unique maximizing
measure for $\phi_0$; and the relatively maximal measures over
$\mu_Y$ are the equilibrium states of $G \circ \pi$,
which can be thought of as the maximizing measures for
$\phi_\infty$.

\section{Bounding the number of ergodic relatively maximal measures}

Let $\pi: X \to Y$ be a 1-block factor map from a 1-step SFT $X$ to a
subshift $Y$ and let $\nu$ be an
ergodic invariant measure on $Y$.
Let $\mu_1 , \dots ,\mu_n \in \cm (X)$ with $\pi \mu_i=
\nu$ for all $i$. Recall the definition of the {\em relatively
  independent joining} $\hat\mu = \mu_1 \otimes \cdots \otimes_\nu
\mu_n$ 
of $\mu_1, \dots , \mu_n$ over $\nu$:
if $A_1, \dots , A_n$ are measurable subsets of $X$ and $\cf$ is the
$\sigma$-algebra of $Y$, then
 \begin{equation}
    \label{eq:relind}
    \hat\mu(A_1\times\ldots\times A_n)=
    \int_Y \prod_{i=1}^n
    \E_{\mu_i}(\one _{A_i}|\pi^{-1}\F)\circ\pi^{-1}
    \,d\nu .
  \end{equation}

  Writing $p_i$ for the projection $X^n\to X$ onto the $i$'th
  coordinate, we note that for $\hat\mu$-almost every $\hat x$ in
  $X^n$, $\pi(p_i(\hat x))$ is independent of $i$. 

We define a number of $\sigma$-algebras on $X^n$.
Denoting by $\cb_X$ the $\sigma$-algebra of $X$ and by $\cb_Y$ the
$\sigma$-algebra of $Y$, let $\cb_0 = \phi^{-1}\cb_Y$, $\cb_i =
p_i^{-1}\cb_X$ for $i=1,\dots ,n$, $\cb_X^-$ the $\sigma$-algebra
generated by $x_n, n<0$, and $\cb_i^-=p_i^{-1}\cb_X^-$ for each $i$. 
{\em Note}: later  we will use the same symbols for corresponding sub-$\sigma$-algebras of 
a different space, $Z = X \times X \times R$.

\begin{defn*}
We say that two measures $\mu_1 , \mu_2 \in \ce (X)$ with $\pi \mu_1 =
\pi \mu_2 = \nu$ are {\em relatively orthogonal} (over $\nu$) and write
  $\mu_1 \perp_\nu \mu_2$ if 
\begin{equation}
(\mu_1 \otimes_\nu  \mu_2) \{ (u,v) \in X \times X : u_0 =v_0\}=0.
\end{equation}
\end{defn*} 

\begin{thm}\label{relorth}
For each ergodic $\nu$ on $Y$, any two distinct ergodic measures on $X$ of
maximal entropy in the fiber $\pi^{-1}\{\nu\}$ are relatively
orthogonal.
\end{thm}

Since $\pi$ is a 1-block
factor map, for each symbol $b$ in the alphabet of $Y$, 
$\pi^{-1}[b]$ consists of a union of 1-block cylinder
sets in $X$. Let $N_\nu(\pi)$ denote the minimum number of cylinders in the
union as $b$ runs over the symbols in the alphabet of $Y$ for which
$\nu [b]>0$.

\begin{cor}\label{bound}
Let $X$ be a $1$-step SFT, $Y$ a subshift on a finite alphabet, and $\pi
: X \to Y$ a $1$-block factor map. 
For any ergodic $\nu$ on $Y$,   
the number of ergodic
  invariant measures of maximal entropy in the fiber
  $\pi^{-1}\{\nu\}$ is at most $N_\nu(\pi)$.
\end{cor}

\begin{proof}
  Suppose that we have $n > N_\nu(\pi)$ ergodic measures $\mu_1,\ldots,\mu_n$
  on $X$, each projecting to $\nu$ and each of maximal entropy in the
  fiber $\pi^{-1}\{ \nu\}$.
Form the relatively independent joining $\hat\mu$ on $X^n$  of the
  measures $\mu_i$ as above.
  Let $b$ be a symbol in the alphabet of $Y$ such that $b$ has
  $N_\nu(\pi)$ preimages $a_1,\dots ,a_{N_\nu(\pi)}$ under the block map $\pi$. 
Since $n > N_\nu(\pi )$, for every $\hat x \in \phi^{-1}[b]$ there are $i
\neq j$ with $(p_i \hat x)_0 = (p_j \hat x)_0$. 
At least one of the sets $S_{i,j} = \{ \hat x \in X^n : (p_i \hat x)_0
= (p_j \hat x)_0 \}$ must have positive $\hat\mu$-measure, and then
also  
 $(\mu_i \otimes _\nu \mu_{j}) \{(u,v) \in X \times X: \pi
u = \pi v, u_0 = v_0 \} > 0$, contradicting Theorem \ref{relorth}. 
\end{proof}

\begin{cor}
Suppose that $\pi: X \to Y$ has a singleton clump: there is a symbol
$a$ of $Y$ whose inverse image is a singleton, which we also denote by
$a$. Then every ergodic measure on $Y$ which 
assigns positive measure to $[a]$ is $\pi$-determinate.
\end{cor}

Before giving the proof of Theorem \ref{relorth}, we 
recall some facts about conditional independence of
$\sigma$-algebras (see \cite[p. 17]{Loeve}) and prove a key lemma.
\begin{lem}\label{condind}
Let $(X,\cb , \mu)$ be a probability space. 
For sub-$\sigma$-algebras $\cb_0, \cb_1, \cb_2$ of $\cb$, the following are
equivalent:
\begin{enumerate}
\item $\cb_1 \perp_{\cb_0} \cb_2$, which is defined by the condition
  that 
for every $\cb_1$-measurable $f_1$ and $\cb_2$-measurable $f_2$,
  $\E (f_1 f_2 | \cb_0) = \E(f_1 | \cb_0) \E (f_2 | \cb_0)$;
\item for every $\cb_2$-measurable $f_2$, $\E (f_2 | \cb_1 \vee \cb_0)
  = \E (f_2 | \cb_0)$;
\item for every $\cb_1$-measurable $f_1$, $\E (f_1 | \cb_2 \vee \cb_0)
  = \E (f_1 | \cb_0)$.
\end{enumerate}
\end{lem}

\begin{lem}\label{condind2}
Let $(X,\cb , \mu)$ be a probability space and let
$\cb_1,\cb_2,\cc_1,\cc_2$ be sub-$\sigma$-algebras of $\cb$.  
If  $\cb_1 \perp_{\cb_0} \cb_2, \cc_1 \subset \cb_1, \cc_2 \subset
\cb_2$, then for every $\cb_1$-measurable $f_1$,
\begin{equation}
\E (f_1 | \cb_0 \vee \cc_1 \vee \cc_2 ) = \E (f_1 | \cb_0 \vee \cc_1 )
.
\end{equation}
\end{lem}
\begin{proof}
First note that $\cb_1 \perp_{\cb_0 \vee \cc_2} \cb_2$, since
for $\cb_1$-measurable $f_1$ we have $\E (f_1 | (\cb_0 \vee \cc_2)
\vee \cb_2 ) = \E (f_1 | \cb_0 \vee \cb_2 ) = \E (f_1 | \cb_0 ) = \E
(f_1 | \cb_0 \vee \cc_2)$. 
Similarly, $\cb_1 \perp_{\cb_0 \vee \cc_1} \cb_2$ and $ \cb_1
\perp_{\cb_0 \vee \cc_1} \cc_2$. 
Thus for any $f_1$ that is $\cb_1$-measurable, $\E ( f_1 | (\cb_0 \vee
\cc_1) \vee \cc_2 ) = \E (f_1 | \cb_0 \vee \cc_1 )$.
\end{proof}

\begin{lem}\label{equa}
Let $\pi : X \to Y$ be a $1$-block factor map from a $1$-step SFT $X$
to a subshift $Y$. Let $\nu$ be an ergodic measure on $Y$ and let
$\mu_1$ and $\mu_2$ be ergodic members of $\pi^{-1}\{ \nu \}$. Let 
$\hat\mu$ be their relatively
independent joining.
If 
$S=\{ (u,v)\in X \times X$ : $u_{-1}=v_{-1}\}$ 
has positive measure with respect to $\hat\mu$ and
for every symbol $j$ in the
alphabet of $X$ 
\begin{equation}
\E_{\hat\mu}(1_{[j]} \circ p_1 | \cb_1^- \vee \cb_0) =
\E_{\hat\mu}(1_{[j]} \circ p_2 | \cb_2^- \vee \cb_0)  \quad\text{a.e.
  on } S ,
\end{equation}
then $\mu_1 = \mu_2$.
\end{lem}
\begin{proof}
Write  $[i]_k$ for the set of points in $X$ whose $k$'th symbol is 
$i$ and $[i]^{(j)}_k$ for $p_j^{-1}[i]_k$. Write $1_{[i]^{(j)}_k}$
for the indicator function of this set. Define
$g^{(j)}_i=\E(1_{[i]^{(j)}_0}|\B_0\vee\B_j^-)$ and set $s_k=\sum_i
1_{[i]^{(1)}_k}1_{[i]^{(2)}_k}=1_{\{ (u,v): u_k=v_k\} }$.
Note that $s_{-1}=1_S$.

Let $\cp$ denote the time-0 partition of $X$ into 1-block cylinder
sets, $\cp_i=p_i^{-1}\cp$ ($i=1,2$) the corresponding partitions
of $X \times X$, and $T=\sigma \times \sigma$.

By assumption, we have $s_{-1}g^{(1)}_i=s_{-1}g^{(2)}_i$ for all
symbols $i$ in the alphabet of $X$.
Taking expectations with respect to $\B_1\vee T\PP_2$, since 
 $s_{-1}g_i^{(1)}$ is $\cb_1 \vee T\PP_2$-measurable, 
we see that 
{\allowdisplaybreaks
\begin{equation}
\begin{aligned}
  s_{-1}g_i^{(1)} &=
  s_{-1}\E(g_i^{(2)}|\B_1\vee T\PP_2)\\
  &=s_{-1}\sum_j \frac{\E(g_i^{(2)}1_{[j]^{(2)}_{-1}}|\B_1)}
  {\E(1_{[j]^{(2)}_{-1}}|\B_1)} 1_{[j]^{(2)}_{-1}}\\
  &=s_{-1}\sum_j \frac{\E(g_i^{(2)}1_{[j]^{(2)}_{-1}}|\B_0)}
  {\E(1_{[j]^{(2)}_{-1}}|\B_0)} 1_{[j]^{(2)}_{-1}},
\end{aligned}
\end{equation}
}
where the last equality follows from Lemma \ref{condind}, noting that
$\cb_0 \subset \cb_1$. Observe that the terms in the
final expression are all measurable with respect to $\B_0\vee T\PP_1
\vee T\PP_2$.

It then follows that
\begin{equation}
  s_{-1}g_i^{(1)}=\E(s_{-1}g_i^{(1)}|\B_0\vee T\PP_1
  \vee T\PP_2)=s_{-1}\E(g_i^{(1)}|\B_0\vee T\PP_1\vee T\PP_2).
\end{equation}
Since $g_i^{(1)}$ is $\B_1$-measurable and $\cb_1$ and $\cb_2$ are
relatively independent over $\cb_0$, by Lemma \ref{condind2} the right
side is equal to $s_{-1}\E(g_i^{(1)}|\B_0\vee T\PP_1)$.
We have thus established the equation
\begin{equation}
  \label{eq:useful}
  s_{-1}\E(g_i^{(1)}|\B_0\vee T\PP_1)=s_{-1}g^{(1)}_i=s_{-1}g^{(2)}_i=
  s_{-1}\E(g_i^{(2)}|\B_0\vee T\PP_2).
\end{equation}

Starting from the equation
$s_{-1}g_i^{(1)}=s_{-1}\E(g_i^{(2)}|\B_0\vee T\PP_2)$, we take
conditional expectations with respect to $\B_1$ to get
\begin{equation}
  \label{eq:gi1}
  \E(s_{-1}|\B_1)g_i^{(1)}=\E(s_{-1}\E(g_i^{(2)}|\B_0\vee
  T\PP_2)|\B_1).
\end{equation}
We have 
\begin{equation}
  \E(g_i^{(2)}|\B_0\vee T\PP_2)=
  \sum_k \frac{\E(g_i^{(2)}1_{[k]^{(2)}_{-1}}|\B_0)}
  {\E(1_{[k]^{(2)}_{-1}}|\B_0)}1_{[k]^{(2)}_{-1}}.
\end{equation}
Hence
\begin{equation}
  s_{-1}\E(g_i^{(2)}|\B_0\vee T\PP_2)=
  \sum_k \frac{\E(g_i^{(2)}1_{[k]^{(2)}_{-1}}|\B_0)}
  {\E(1_{[k]^{(2)}_{-1}}|\B_0)}1_{[k]^{(1)}_{-1}}1_{[k]^{(2)}_{-1}}.
\end{equation}

Substituting this in (\ref{eq:gi1}) and again using relative
independence, we see that
{\allowdisplaybreaks 
\begin{equation} 
\begin{aligned}
  \E(s_{-1}|\B_1)g_i^{(1)}&=
  \sum_k \frac{\E(g_i^{(2)}1_{[k]^{(2)}_{-1}}|\B_0)}
  {\E(1_{[k]^{(2)}_{-1}}|\B_0)}1_{[k]^{(1)}_{-1}}\E(1_{[k]^{(2)}_{-1}}|\B_1)\\
  &=\sum_k \E(g_i^{(2)}1_{[k]^{(2)}_{-1}}|\B_0)1_{[k]^{(1)}_{-1}}.
\end{aligned}
\end{equation}
}
We observe that the right-hand side and also $\E(s_{-1}|\B_1)$ are
$\B_0\vee T\PP_1$-measurable (using the definition of $s_{-1}$ and
relative independence).
Hence provided that $\E(s_{-1}|\B_1)>0$ a.e., we will have that  $g_i^{(1)}$
is $\B_0\vee T\PP_1$-measurable, and similarly  $g_i^{(2)}$
is $\B_0\vee T\PP_2$-measurable.

We now demonstrate that $\E(s_{-1}|\B_1)>0$ on a set of full measure.
To prove this, we note that $\E(s_{-1}|\B_1)$ is of the form $f\circ
p_1$ for $f$ a function on $X$. Thus if we can show that
$\E(s_{-1}|\B_1)(x)>0$ implies $\E(s_{-1}|\B_1)(Tx)>0$, it will follow
that the set where $f$ is positive is invariant and hence of measure 0
or 1 by ergodicity of $\mu_1$. Since the integral of the function is
positive (being equal to $(\mu_1 \otimes_\nu \mu_2 )\{ (u,v): u_{-1} =
v_{-1} \}$), to show 
that the function is positive on a set of full measure it is enough to
establish the above invariance.

Now
{\allowdisplaybreaks
\begin{equation}\begin{aligned}
\E(s_{-1}|\B_1)(Tx)&=\E(s_0|\B_1)(x) \\
 &=\sum_{i}\E(1_{[i]^{(1)}_0}1_{[i]^{(2)}_0}|\B_1)\\
  &=\sum_i1_{[i]^{(1)}_0}\E(1_{[i]^{(2)}_0}|\B_1)\\
  &\geq\sum_i1_{[i]^{(1)}_0}\E(s_{-1}1_{[i]^{(2)}_0}|\B_1)\\
  &=\sum_i1_{[i]^{(1)}_0}\E(\E(s_{-1}1_{[i]^{(2)}_0}|\B_1\vee
  T\PP_2)|\B_1).
\end{aligned}
\end{equation}
Using Lemma \ref{condind2}, this equals
\begin{equation}
\begin{aligned}
  &\sum_i1_{[i]^{(1)}_0}\E(s_{-1}\E(1_{[i]^{(2)}_0}|\B_1\vee
  T\PP_2)|\B_1)\\ 
  &=\sum_i1_{[i]^{(1)}_0}\E(s_{-1}\E(1_{[i]^{(2)}_0}|\B_0\vee
  T\PP_2)|\B_1)\\
  &=\sum_i1_{[i]^{(1)}_0}\E(s_{-1}g_i^{(2)}|\B_1)  \text{  (from (\ref{eq:useful}))}\\
  &=\sum_i1_{[i]^{(1)}_0}\E(s_{-1}g_i^{(1)}|\B_1)\\
  &=\sum_i g_i^{(1)}1_{[i]^{(1)}_0}\E(s_{-1}|\B_1)\\
  &=\E(s_{-1}|\B_1)\sum_i 1_{[i]_0^{(1)}}\E(1_{[i]_0^{(1)}}|\B_0 \vee
  \B_1^-).
\end{aligned}
\end{equation}
}
For $x$ in a set of full measure, 
$1_D(x)>0$ implies ${\E}(1_D|{\mathcal{F}})(x)>0$
(consider integrating the conditional expectation over the set where
it takes the value 0), so the sum on the right-hand side
of the above is positive almost everywhere. Since the first factor is
positive by assumption, the conclusion that $\E(s_0|\B_1)>0$ follows,
allowing us to deduce that $g_i^{(j)}$ is $\B_0\vee T\PP_j$-measurable.

Now we may write $g_i^{(j)}$ as
\begin{equation}
  g_i^{(j)}=\sum_k 1_{[k]^{(j)}_{-1}} h_{k,i}^{(j)},
\end{equation}
where the $h_{k,i}^{(j)}$ are $\cb_0$-measurable.
Writing out the equation $s_{-1}g^{(1)}_i=s_{-1}g^{(2)}_i$, we have
\begin{equation}
  \sum_k 1_{[k]^{(1)}_{-1}}1_{[k]^{(2)}_{-1}} h_{k,i}^{(1)}=
  \sum_k 1_{[k]^{(1)}_{-1}}1_{[k]^{(2)}_{-1}} h_{k,i}^{(2)}.
\end{equation}

Since for distinct $k$, the terms are disjointly supported, we have
for each $k$, 
\begin{equation}
  1_{[k]^{(1)}_{-1}}1_{[k]^{(2)}_{-1}} h_{k,i}^{(1)}=
  1_{[k]^{(1)}_{-1}}1_{[k]^{(2)}_{-1}} h_{k,i}^{(2)}.
\end{equation}
Taking conditional expectations of both sides with respect to $\B_0$
and using Lemma \ref{condind},
we deduce
\begin{equation}
  \E(1_{[k]^{(1)}_{-1}}|\B_0)\E(1_{[k]^{(2)}_{-1}}|\B_0)
  (h_{k,i}^{(1)}-h_{k,i}^{(2)})=0\quad\text{a.e.}
\end{equation}
From this we see that if $\E(1_{[k]^{(1)}_{-1}}|\B_0)>0$ and
$\E(1_{[k]^{(2)}_{-1}}|\B_0)>0$, then $h_{k,i}^{(1)}=h_{k,i}^{(2)}$.
This allows us to make the following definition:
\begin{equation}
  h_{k,i}=
  \begin{cases}
    h_{k,i}^{(1)}&\text{if $\E(1_{[k]^{(1)}_{-1}}|\B_0)>0$}\\
    h_{k,i}^{(2)}&\text{if $\E(1_{[k]^{(2)}_{-1}}|\B_0)>0$ .}
  \end{cases}
\end{equation}
It follows that
\begin{equation}
  g_i^{(j)}=\sum_k h_{k,i}1_{[k]^{(j)}_{-1}} \quad\hat\mu\text{-a.e..}
\end{equation}

We now show that the two measures agree.
We will show by induction on the length of the cylinder set that for
any $\B_0$-measurable function $f$ and any cylinder set $C$ in $X$,
\begin{equation}
\int 1_S 1_C\circ p_1 f \,d\hat\mu=
\int 1_S 1_C\circ p_2 f \,d\hat\mu.
\end{equation}
To start the induction, let $C$ be the cylinder set $[i_0]$ in
$X$. Then
\begin{equation}
\begin{aligned}
  \int 1_S 1_{[i_0]^{(j)}} f\,d\hat\mu&=
  \int 1_S f \E(1_{[i_0]^{(j)}}|\B_0\vee \B_1^- \vee \B_2^-)\,d\hat\mu\\
  &=\int 1_S f g_{i_0}^{(j)}\,d\hat\mu  ;
\end{aligned}
\end{equation}
but by assumption $1_S g_i^{(1)}=1_S g_i^{(2)}$, 
showing the result in the case that $C$ is a cylinder of length 1.
Now suppose that the result holds for cylinders of length $n$ and let
$C=[i_0\ldots i_n]$. Write $D=[i_0\ldots i_{n-1}]$.
Now
\begin{equation}
\begin{aligned}
  \int 1_S (1_C\circ p_j) f\,d\hat\mu&=
  \int 1_S (1_D\circ p_j) 1_{[i_n]^{(j)}_n} f\,d\hat\mu\\
  &=\int 1_S (1_D\circ p_j) f\E(1_{[i_n]^{(j)}_n} |T^{-n}\B_1^-\vee
  T^{-n}\B_2^-\vee\B_0)\,d\hat\mu\\
  &=\int 1_S (1_D \circ p_j) f g_{i_n}^{(j)}\circ T^n\,d\hat\mu\\
  &=\int 1_S (1_D \circ p_j) f h_{i_{n-1},i_n}\circ T^n\,d\hat\mu.
\end{aligned}
\end{equation}
Since $h_{i_{n-1},i_n}\circ T^n$ is $\B_0$-measurable, it follows from
the induction hypothesis that the integrals are equal for $j=1$ and $j=2$
as required.

In particular, taking $f$ to be 1, we have
$\hat\mu(S\cap p_1^{-1}C)=\hat\mu(S\cap p_2^{-1}C)$ for all $C$.
Letting $\hat\nu(A)=\hat\mu(S\cap A)$, we see that
$\hat\nu\circ p_1^{-1}=\hat\nu\circ p_2^{-1}$.
Since $\mu_i(A)\geq \hat\nu\circ p_i^{-1}(A)$ for all $A$ and the
measures $\mu_i$ are ergodic, it follows that $\mu_1$ and $\mu_2$ are
not mutually singular and hence are equal.
\end{proof}


\begin{proof}[Proof of Theorem \ref{relorth}]
Let $\mu_1$ and $\mu_2$ be two different ergodic relatively maximal
measures over $\nu \in \ce (Y)$ 
and suppose that they are not relatively orthogonal, 
so that $(\mu_1 \otimes_\nu  \mu_2) \{ (u,v) \in X \times X : u_0 =v_0\}>0$.
Let $\hat\mu = \mu_1 \otimes_\nu \mu_2$.
We will construct a measure on $X$ with strictly greater entropy than
  $\mu_1$ or $\mu_2$ by building a larger space from which the new
  measure will appear as a factor. 
 (J. Steif reminded us that a similar interleaving of two processes is used in
  \cite{FurstenbergPW} for a different purpose.)
  
  Let $R$ denote the set $\{1,2\}^{\Z}$, and let $\beta$ be
  the Bernoulli measure on $R$ with probabilities
  $\frac{1}{2},\frac{1}{2}$. Write $(r_n)_{n\in\Z}$ for a
  typical element of $R$.
  Form $Z=X^2\times R$ with invariant
  measure $\eta=\ri\times\beta$. We then define maps from $Z$ to $X$ as
  follows. 
  Given a point $(u,v,r)\in Z$, set $\pi_1(u,v,r)=u$, $\pi_2(u,v,r)=v$
  and write $N_k(u,v)$ for $\sup\{n<k\colon u_n=v_n\}$. Note that this
  quantity may be $-\infty$ if there are no coincidences. We will take
  $r_{-\infty}$ to be a further random variable taking the values 1
  and 2 with equal probability for each $r\in R$.
  Define $\pi_3\colon Z\to X$ by
  \begin{equation}
    \pi_3(u,v,r)_k=
    \begin{cases}
      u_k&\text{if $r_{N_k(u,v)}=1$}\\
      v_k&\text{if $r_{N_k(u,v)}=2$ .}
    \end{cases}
  \end{equation}

  To see that $\pi_3(u,v,r)$ is indeed a point of $X$, note that it
  consists of concatenations of parts of $u$ and $v$, changing only at
  places where they agree. As a corollary, since $\pi(u)=\pi(v)$ for
  almost all $(u,v,r)\in Z$, it follows that
  $\pi(\pi_3(z))=\pi(\pi_2(z))=\pi(\pi_1(z))$ for $\eta$-almost
  every $z$ in $Z$. Write $\Phi$ for the factor mapping $\pi\circ
  \pi_1$ from $(Z,\eta)$ to $(Y,\nu)$.

  By construction $\mu_1=\eta\circ\pi_1^{-1}$ and
  $\mu_2=\eta\circ\pi_2^{-1}$. Define $\mu_3=\eta\circ\pi_3^{-1}$. We
  shall then demonstrate that
  $h_{\mu_3}(X)>h_{\mu_1}(X)=h_{\mu_2}(X)$.

  We define $\sigma$-algebras on $Z$ corresponding to those appearing above. Letting
  $\B_X$ be the Borel $\sigma$-algebra on $X$ as before, we set for
  each $i=1,2,3$, $\B_i=\pi_i^{-1}\B_X$. Write
  $\B_X^-$ for the $\sigma$-algebra generated by the cylinder
  sets in $X$ depending on coordinates $x_n$ for $n<0$. These then
  give $\sigma$-algebras $\B_i^-$ on $Z$ defined by
  $\B_i^-=\pi_i^{-1}\B_X^-$. We will require two
  further $\sigma$-algebras, $\B_0=\Phi^{-1}\B_Y$
  with $\B_0^-$ being defined analogously to the above.
Note that $\B_i \supset \B_0$ for $i=1,2,3$. 
  
  Again reusing previous notation in a slightly different context,
 continue to denote by $\P$ the partition of $X$
  into time 0 cylinders and write $\P_i$ for
  $\pi_i^{-1}\P$, so that for $i=1,2,3$, $\P_i$ is a
  partition of $Z$. Finally, write
  $\Q=\Phi^{-1}\{[j]\colon\text{$[j]$ is a cylinder set in
  $Y$}\}$.

  It is useful to note the following property of (\ref{eq:relind}):
 If $A_1\in\B_1$ and $A_2\in\B_2$,
  then 
  \begin{equation}
    \label{eq:relind2}
    \eta(A_1\cap A_2)=\int
    \E_\eta(1_{A_1}|\B_0)
    \E_\eta(1_{A_2}|\B_0)
    \,d\eta.
  \end{equation}

 We will use the fact that
if $f$ is $\B_1$-measurable, then
  \begin{equation}
    \label{eq:key}
    \E_\eta(f|\B_2)=\E_\eta(f|\B_0),
  \end{equation}
a consequence of Lemma \ref{condind}.

  Standard results of entropy theory tell us that
  $h_{\mu_i}(X)=H_\eta(\P_i|\B_i^-)$. Further, by
  Pinsker's Formula (see \cite[theorem 6.3, p. 67]{Parry69}, applied with $\beta$
  coarser than $\alpha$), this can be
  re-expressed as 
  \begin{equation}
    \label{eq:ent}
    h_{\mu_i}(X)=H_\eta(\P_i|
    \B_i^-\vee\B_0)+H_\eta(\Q|\B_0^-)=H_\eta(\P_i|\B_i^-\vee
    \B_0)+h_\nu(Y).
  \end{equation}
  Since $\mu_1$ and $\mu_2$ were presumed
  to be measures of maximal entropy in the fiber, they have equal
  entropy and hence $H_\eta(\P_1|
  \B_1^-\vee\B_0)=H_\eta(\P_2|
  \B_2^-\vee\B_0)$.  Our aim is to show that this
  leads to a contradiction by showing that
  $H_\eta(\P_3|
  \B_3^-\vee\B_0)>H_\eta(\P_1|\B_1^-\vee
  \B_0)$.  By definition,
  \begin{equation}
    \label{eq:entropy}
    \begin{gathered}
      H_\eta(\P_i|\B_i^-\vee\B_0)=
\int -\sum_{j}(\one _{[j]}\circ\pi_i)
      \log \E(\one _{[j]}\circ\pi_i| \B_i^-\vee B_0)\,d\eta\\
      =\int -\sum_{j}\E(\one _{[j]}\circ\pi_i|\B_i^-\vee B_0)
      \log \E(\one _{[j]}\circ\pi_i|\B_i^-\vee B_0)\,d\eta\\
      =\int \sum_j \psi\left(\E(\one _{[j]}\circ\pi_i|\B_i^-\vee
      B_0)\right)\,d\eta,
    \end{gathered}
  \end{equation}
  where $\psi$ is the strictly concave function $[0,1]\to[0,1]$,
  $\psi(x)=-x\log x$ (with $\psi(0)$ defined to be 0).

    The following claim is an essential point of the argument.
  We shall show that
{\allowdisplaybreaks
  \begin{equation}
    \label{eq:g3}
    \begin{gathered}
      \E_\eta(\one _{[j]}\circ\pi_3|\B_1^-\vee\B_2^-\vee\B_3^-\vee\B_0)(z)=\\
              \E_\eta(\one _{[j]}\circ\pi_1|\B_1^-\vee\B_0)
        \quad\text{if $\pi_3(z)_{-1}=\pi_1(z)_{-1}\neq\pi_2(z)_{-1}$;}\\
        \E_\eta(\one _{[j]}\circ\pi_2|\B_2^-\vee\B_0)
        \quad\text{$\pi_3(z)_{-1}=\pi_2(z)_{-1}\neq\pi_1(z)_{-1}$;}\\
        \frac{1}{2}\E_\eta(\one _{[j]}\circ\pi_1|\B_1^-\vee\B_0)+
        \frac{1}{2}\E_\eta(\one _{[j]}\circ\pi_2|\B_2^-\vee\B_0)
        \quad\text{$\pi_3(z)_{-1}=\pi_1(z)_{-1}=\pi_2(z)_{-1}$.}
          \end{gathered}
  \end{equation}
}
  Clearly, the right-hand side of the equation is measurable with
  respect to $\B_1^-\vee\B_2^-\vee\B_3^-\vee\B_0$. To verify the
  claim, it will be sufficient to integrate the right-hand side over
  the elements of a
  generating semi-algebra of $\B_1^-\vee\B_2^-\vee\B_3^-\vee\B_0$.
  Specifically, we will integrate over sets of the form $A\cap B\cap
  C\cap D$, where $A$, $B$ and $C$ are the preimages under the
  respective maps of cylinder sets in $X$ of a common length (ending
  at time $-1$) and $D \in \B_0$. 

  Suppose $A$, $B$, and $C$ are cylinders depending on
  the coordinates $-n$ to $-1$ of $\pi_1(z)$, $\pi_2(z)$, and
  $\pi_3(z)$ and 
  that $A\cap B\cap C$ has
  positive measure.
 Then for $z\in A\cap B\cap C$, $\pi_3(z)_{-1}$ is
  equal to either $\pi_1(z)_{-1}$ or $\pi_2(z)_{-1}$ (or both) by
  definition of $\pi_3$. Further, $\pi_1(z)_{-1}$, $\pi_2(z)_{-1}$, and
  $\pi_3(z)_{-1}$
 are constant over 
  the intersection in question.

  If on $A\cap B\cap C$,
  $\pi_3(z)_{-1}=\pi_1(z)_{-1}\neq\pi_2(z)_{-1}$,
  then we calculate
  \begin{equation}\begin{aligned}
    &\int_{A\cap B\cap C\cap D}
    \E_\eta(\one _{[j]}\circ\pi_1|\B_1^-\vee \B_0)(z)\,d\eta=\\
    &\int \one _B \one _C \E_\eta(\one _{[j]}\circ\pi_1\one _A\one _D|\B_1^-\vee
    \B_0)\,d(\ri\times\beta).
  \end{aligned}\end{equation}
  Performing first the integration over $R$ with respect to the
  measure $\beta$, we see that the only factor depending on the random
  part $r\in R$ is $1_C$, the others being functions only of $(u,v)
  \in X^2$. The coordinates of
  $\pi_3(z)$ from $-n$ to $-1$ are concatenations of blocks of
  $\pi_1(z)$ and $\pi_2(z)$, the choice (between a block in $u$ and a
  {\em different} block in $v$) being made according to the
  entries in $r$, hence with probabilities $1/2, 1/2$). 
 If $k=k_{A,B}(u,v) = 1 + \text{card}\{j: -n \leq j \leq -2,
u_j=v_j,  u_{j+1} \neq v_{j+1}\}$, then
\begin{equation}
\int_R \one_C(u,v,r) \, d\beta (r) = \frac{1}{2^{k_{A,B}(u,v)}},
\end{equation}
which is constant on $A \cap B$. 
The following calculation will be more readable if we write $\E^\cb f$
for $\E (f|\cb )$.
Since $B \in \cb_2$ and $\cb_2 \perp_{\cb_0}\cb_1^-$, we have
$\E^{\cb_1^- \vee \cb_0}\one_B = \E^{\cb_0}\one_B$.
Consequently,
{\allowdisplaybreaks
  \begin{equation}\begin{aligned}
    &\int_{A\cap B\cap C\cap D}
   \E_\eta(\one _{[j]}\circ\pi_1|\B_1^-\vee \B_0)(z)\, d\eta\\
&=  \int_{X^2} 2^{-k}\one_D \one_B \one_A \E^{\cb_1^- \vee
  \cb_0}(\one_{[j]} \circ \pi_1) \, d\ri(u,v)\\
&= \int_{X^2} 2^{-k}  \one_B \E^{\cb_1^- \vee
  \cb_0}(\one_D \one_A \cdot (\one_{[j]} \circ \pi_1)) \, d\ri(u,v)\\
&= \int_{X^2} 2^{-k}  \E^{\cb_1^- \vee  \cb_0}[
\one_B \E^{\cb_1^- \vee
  \cb_0}(\one_D \one_A \cdot (\one_{[j]} \circ \pi_1))] \, d\ri(u,v)\\
&= \int_{X^2} 2^{-k}  [\E^{\cb_1^- \vee
  \cb_0} \one_B] [\E^{\cb_1^- \vee
  \cb_0}(\one_D \one_A \cdot (\one_{[j]} \circ \pi_1))] \, d\ri(u,v)\\
&= \int_{X^2} 2^{-k}  [\E^{ \cb_0} \one_B] [\E^{\cb_1^- \vee
  \cb_0}(\one_D \one_A \cdot (\one_{[j]} \circ \pi_1))] \, d\ri(u,v)\\
&= \int_{X^2} 2^{-k}  [\E^{ \cb_0} (\one_B \one_D)] [\E^{\cb_1^- \vee
  \cb_0}(\one_A \cdot (\one_{[j]} \circ \pi_1))] \, d\ri(u,v)\\
&= \int_{X^2} 2^{-k}  \E^{\cb_0} \{ [\E^{ \cb_0} (\one_B \one_D)] [\E^{\cb_1^- \vee
  \cb_0}(\one_A \cdot (\one_{[j]} \circ \pi_1))] \} \, d\ri(u,v)\\
&= \int_{X^2} 2^{-k}  [\E^{ \cb_0} (\one_B \one_D)] [(\E^{ \cb_0}
(\one_A \cdot (\one_{[j]} \circ \pi_1))] \, d\ri(u,v)\\
&=\eta(A\cap B\cap C\cap D\cap \pi_1^{-1}[j])
    =\eta(A\cap B\cap C\cap D\cap \pi_3^{-1}[j]),
 \end{aligned}\end{equation}
}
by (\ref{eq:relind2}), since $B,D \in \cb_2$ and $A, \pi_1^{-1}[j]
\in \cb_1$.
  This demonstrates the desired equality in the case
  $\pi_3(z)_{-1}=\pi_1(z)_{-1}\neq \pi_2(z)_{-1}$. The case
  $\pi_3(z)_{-1}=\pi_2(z)_{-1}\neq \pi_1(z)_{-1}$ is dealt with
  similarly.

  If $\pi_3(z)_{-1}=\pi_1(z)_{-1}=\pi_2(z)_{-1}$, then the integrand
  is the average of the two previous integrands, so we see that
  \begin{equation}\begin{aligned}
    &\int_{A\cap B\cap C\cap D}\left(
    \textstyle{\frac{1}{2}}\E_\eta(\one _{[j]}\circ\pi_1|\B_1^-\vee\B_0)+
    \textstyle{\frac{1}{2}}\E_\eta(\one _{[j]}\circ\pi_2|\B_2^-\vee\B_0)
    \right)\,d\eta=\\
    &\textstyle{\frac{1}{2}}\eta(A\cap B\cap C\cap D\cap \pi_1^{-1}[j])+
    \textstyle{\frac{1}{2}}\eta(A\cap B\cap C\cap D\cap
    \pi_2^{-1}[j])=\\
    &\eta(A\cap B\cap C\cap D\cap \pi_3^{-1}[j]).
  \end{aligned}\end{equation}
  This completes the proof of equation (\ref{eq:g3}).

  Using (\ref{eq:entropy}), we have
  \begin{equation}\begin{aligned}
    H_\eta(\P_3|\B_3^-\vee\B_0)&
    \geq H_\eta(\P_3|\B_1^-\vee\B_2^-\vee\B_3^-\vee\B_0)\\
    &=\int\sum_j\psi(
    \E_\eta(\one _{[j]}\circ\pi_3|\B_1^-\vee\B_2^-\vee\B_3^-\vee\B_0))\,d\eta.
  \end{aligned}\end{equation}
  We separate the integral into parts according to whether
  $\pi_3(z)_{-1}$ is equal to $\pi_1(z)_{-1}$, $\pi_2(z)_{-1}$ or
  both. Let $S_1=\{z\colon \pi_3(z)_{-1}=\pi_1(z)_{-1}\neq
  \pi_2(z)_{-1}\}$, $S_2=\{z\colon \pi_3(z)_{-1}=\pi_2(z)_{-1}\neq
  \pi_1(z)_{-1}\}$ and $S_3=\{z\colon\pi_3(z)_{-1}=\pi_1(z)_{-1}=
  \pi_2(z)_{-1}\}$. Let $A=\{z\colon \pi_1(z)_{-1}\neq
  \pi_2(z)_{-1}\}$ so that $A=S_1\cup S_2$. Note that $S_1$ and $S_2$
  have equal measure by definition of $\pi_3$.

By symmetry,
\begin{equation}
\int_{S_1}
    \sum_j\psi(\E_\eta(\one _{[j]}\circ\pi_1|\B_1^-\vee\B_0)\,d\eta =
\int_{S_2}
    \sum_j\psi(\E_\eta(\one _{[j]}\circ\pi_1|\B_1^-\vee\B_0)\,d\eta,
\end{equation}
so by (\ref{eq:g3}),
  \begin{equation}\begin{aligned}
    &\int_{S_1}\sum_j\psi(
    \E_\eta(\one _{[j]}\circ\pi_3|\B_1^-\vee\B_2^-\vee\B_3^-\vee\B_0))\,d\eta=\\
    &\textstyle{\frac{1}{2}}\int_{A}
    \sum_j\psi(\E_\eta(\one _{[j]}\circ\pi_1|\B_1^-\vee\B_0)\,d\eta.
  \end{aligned}\end{equation}
Similarly,
  \begin{equation}\begin{aligned}
    &\int_{S_2}\sum_j\psi(
    \E_\eta(\one _{[j]}\circ\pi_3|\B_1^-\vee\B_2^-\vee\B_3^-\vee\B_0))\,d\eta=\\
    &\textstyle{\frac{1}{2}}\int_{A}
    \sum_j\psi(\E_\eta(\one _{[j]}\circ\pi_2|\B_2^-\vee\B_0)\,d\eta.
  \end{aligned}\end{equation}
  Finally, integrating over $S_3$,
  \begin{equation}\begin{aligned}
    &\int_{S_3}\sum_j\psi(\E_\eta(\one _{[j]}\circ\pi_3|
    \B_1^-\vee\B_2^-\vee\B_3^-\vee\B_0))\,d\eta=\\
    &\int_{A^c}\sum_j\psi\left(
      \textstyle{\frac12}(\E_\eta(\one _{[j]}\circ\pi_1|\B_1^-\vee\B_0)
      +\E_\eta(\one _{[j]}\circ\pi_2|\B_2^-\vee\B_0))\right)\,d\eta>\\
    &\textstyle{\frac12}\int_{A_c}\sum_j\left(
      \psi(\E_\eta(\one _{[j]}\circ\pi_1|\B_1^-\vee\B_0)+
     \psi(\E_\eta(\one _{[j]}\circ\pi_2|\B_2^-\vee\B_0)\right)\,d\eta.
  \end{aligned}\end{equation}

  The strict inequality in the above arises since $\psi$ is
  strictly concave and there exist a $j$ in the alphabet of $X$ and a set of points of positive measure in
  $A^c=\{(u,v,r)\in Z = X^2 \times R: u_{-1}=v_{-1} \}$ for which 
  $\E_\eta(\one _{[j]}\circ\pi_1|\B_1^-\vee\B_0) \neq
  \E_\eta(\one _{[j]}\circ\pi_2|\B_2^-\vee\B_0)$---
for, if not, Lemma \ref{equa} would imply that $\mu_1 = \mu_2$.

  Now adding the preceding equalities, we see
  \begin{equation}\begin{aligned}
    &H_\eta(\P_3|\B_3^-\vee\B_0)>\\
    &\frac{1}{2}
    \left(\int \sum_j\psi(\E_\eta(\one _{[j]}\circ\pi_1|\B_1^-\vee\B_0)\,d\eta+
      \int\sum_j\psi(\E_\eta(\one _{[j]}\circ\pi_2|\B_2^-\vee\B_0)
    \,d\eta\right)\\
    &=\textstyle{\frac12}(H(\P_1|\B_1^-\vee\B_0)+H(\P_2|\B_2^-\vee\B_0))\\
    &=h_{\mu_1}(X)-h_\nu(Y) .
  \end{aligned}\end{equation}
  From (\ref{eq:ent}), we see that $h_{\mu_3}(X)>h_{\mu_1}(X)$ as required.
\end{proof}

\begin{rem*}\label{smb}
It would be desirable to
have  a proof of this result based on 
the Shannon-McMillan-Breiman Theorem, but so far we have not been
able to construct one.
\end{rem*}

\begin{defn*}
Let $(X, \cb, \mu ,T)$ and $(Y, \cc, \nu ,S)$ be measure-preserving
systems, $\pi: X \to Y$ a factor map, and
$\alpha$ a finite
generating partition for $X$. We say that $\mu$ is {\em relatively
  Markov for $\alpha$ over $Y$} if it satisfies one of the following
two equivalent conditions:
\begin{enumerate}
\item $\alpha \perp_{T^{-1}\alpha \vee \pi^{-1}\cc} \alpha_2^\infty$ ;
\item $H_\mu (\alpha | \alpha_1^\infty \vee \pi^{-1}\cc) = H_\mu (\alpha |
  T^{-1}\alpha \vee \pi^{-1}\cc)$.
\end{enumerate}
(As usual, $\alpha_i^j= \bigvee_{k=i}^j T^{-k}\alpha$.)
\end{defn*}

\begin{cor}
If $X$ is a $1$-step SFT, $Y$ is a subshift, $\pi : X \to Y$ is a $1$-block factor map, $\nu$ is an
ergodic measure on $Y$, and $\mu$ is an ergodic relatively maximal
measure over $\nu$, then $\mu$ is relatively Markov for the time-0
partition of $X$ over $Y$.
\end{cor}
\begin{proof}
We apply the first half of the proof of Lemma \ref{equa} with
$\mu_1=\mu_2=\mu$. Note that then $\hat\mu (S) >0$.
If $s_{-1}g_i^{(1)}=s_{-1}g_i^{(2)}$ for all symbols $i$ in the alphabet
of $X$, the proof proceeds as before to show that 
the information function with respect
to $\mu$ of the time-0
partition $\cp$ of $X$ given $\cp_1^\infty \vee \pi^{-1}\cb_Y$ is measurable with
respect to $\cp \vee\sigma^{-1}\cp \vee \pi^{-1} \cb_Y$, 
and hence $\mu$ is a 1-step relatively Markov measure.

If there is a symbol $i$ in the alphabet of $X$ for which 
$s_{-1}g_i^{(1)} \neq s_{-1}g_i^{(2)}$, 
then the construction in the proof of Theorem \ref{relorth}, by
interleaving strings according to another random process, will again
produce a measure projecting to $\nu$ which will have entropy greater
than $h(\mu)$.
\end{proof}

\section{Examples}

\begin{ex}
In case $\pi$ has a singleton clump $a$ and $\nu$ is Markov on $Y$, we can
construct the unique relatively maximal measure above $\nu$
explicitly.  
Denote the cylinder sets $[a]$ in $X$ and in $Y$ by $X_a$ and $Y_a$,
respectively. If $\nu$ is (1-step) Markov on $Y$, then the
first-return map $\sigma_a: Y_a \to Y_a$ is countable-state Bernoulli
with respect to the restricted and normalized measure $\nu_a=\nu
/\nu[a]$: the states are all the loops or return blocks $aC^i$ with 
$aC^ia=ac^i_1\dots c^i_{r_i}a$ appearing in $Y$ and no $c^i_j=a$.

Under $\pi^{-1}$, the return blocks to $[a]$ expand into bands
$aB^{i,j}$, with $aB^{i,j}a$ appearing in $X$ and 
$\pi B^{i,j}=C^i$ for all $i,j$. 
Topologically, $(X_a,\sigma_a )$ is a countable-state full shift on
these symbols $aB^{i,j}$. We define $\mu_a$ to be the countable-state
Bernoulli measure on $(X_a,\sigma_a)$ which equidistributes the measure 
of each loop (state) of $Y_a$ over its preimage band:
\begin{equation}
\mu_a[aB^{i,j}]= \frac{\nu_a[aC^ia]}{|\pi^{-1}[aC^ia]|} \quad\text{for
  all } i,j .
\end{equation}

We show now that this choice of $\mu_a$ is relatively maximal over
$\nu_a$. Let $\lambda_a$ be any probability measure on $X_a$ which
maps under $\pi$ to $\nu_a$. Then the countable-state Bernoulli
measure on $X_a$ which agrees with $\lambda_a$ on all the 1-blocks
$aB^{i,j}$ (its ``Bernoullization'') 
has entropy no less than that of $\lambda_a$ and still
projects to the Bernoulli measure $\nu_a$, so we may as well assume
that $\lambda_a$ is countable-state Bernoulli.
If $\lambda_a[aB^{i,j}]=q^{i,j}$ and $|\pi^{-1}(aC^ia)|=J_i$ for all
$i,j$, then
\begin{equation}
h(X_a,\sigma_a,\lambda_a) = \sum_{i=1}^\infty \sum_{j=1}^{J_i}
q^{i,j}\log q^{i,j} .
\end{equation}
Note that for each $i$
\begin{equation}
\sum_{j=1}^{J_i} q^{i,j} = \nu_a[aC^ia]
\end{equation}
is fixed at the same value for all $\lambda_a$.
Thus for each $i$,
\begin{equation}
\sum_{j=1}^{J_i} q^{i,j}\log q^{i,j}
\end{equation}
is maximized by putting all the $q^{i,j}$ equal to one another.

Finally, this unique relatively maximal $\mu_a$ over $\nu_a$ determines
the unique relatively maximal $\mu$ on $X$ over $\nu$ on $Y$, since
according to 
Abramov's formula
\begin{equation}
h(X,\sigma ,\mu ) = \mu[a]\, h(X_a,\sigma_a,\mu_a),
\end{equation}
and $\mu[a]=\nu[a]$.

We show how this calculation of the unique relatively maximal measure
over a Markov measure in the case of a singleton clump works out in a
particular case. 
It was shown in \cite{Shin, Shin2} that for the following factor map there is a 
saturated compensation function $G \circ \pi$ with $G \in \cc (Y)$
but no such compensation function with $G \in \cf (Y)$. There is a
singleton clump, $a$.

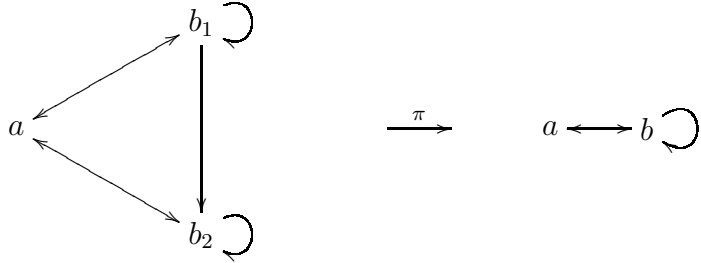
\bigskip
\begin{center}
\xymatrix{
& & & &     b_1 \ar@(ur,dr)[] \ar[dd]  \\
& & a\ar@{<->}[urr] \ar@{<->}[drr] & & & & \ar [r]^\pi & & a \ar@{<->}
[r] & b \ar@(ur,dr) \\
& & & &    b_2 \ar@(ur,dr)[]            }
\end{center}

\bigskip
\noindent
For each $k \geq 1$ the block $ab^ka$ in $Y$ has 
$k+1$ preimages, depending on when the subscript on $b$ switches from
1 to 2. Let $\nu$ be Markov on $Y$. 
To each preimage $aB_1aB_2a\dots aB_r$ of
$ab^{k_1}ab^{k_2}\dots ab^{k_r}$ the optimal measure $\mu_a$ assigns
  measure
\begin{equation}
\mu_a[aB_1aB_2a\dots aB_r] = \frac{1}{k_1+1} \dots
\frac{1}{k_r+1}\nu_a[ ab^{k_1}ab^{k_2}\dots ab^{k_r}].
\end{equation}
The unique relatively maximal
measure over $\nu_a$ can be described in terms of fiber measures as
follows. Given $y=ab^{k_1}ab^{k_2}\dots ab^{k_r}\dots \in Y_a$, $\mu_{a,y}$
chooses the preimages of each $b^{k_i}$ with equal probabilities and
independently of the choice of preimage of any other $b^{k_j}$. Then
\begin{equation}
\mu_a[aB_1aB_2a\dots aB_r] = \int_{Y_a} \mu_{a,y}[aB_1aB_2a\dots aB_r]\, d\nu_a
(y) .
\end{equation}
\end{ex}

\begin{ex}
The relatively maximal measures over 
an ergodic measure $\nu$ on $Y$ which is supported on the orbit 
$\mathcal O (y)$ of a periodic point  $y=CCC\dots \in Y$
can be found by analyzing the SFT $X_y=\pi^{-1}\mathcal O(y)$. 
The relatively
maximal measures over $\nu$ are determined by the maximal
(Shannon-Parry) measures on the irreducible components of $X_y$. 
Consequently, if $X_y$ is irreducible, then the
discrete invariant measure on the orbit of $y$ is $\pi$-determinate.
\end{ex}

\begin{ex}\label{Walters-nopidet}
{\em Failure of $\pi$-determinacy for a fully-supported measure}.
In the preceding example, along with others discussed in
\cite{Petersen00},  failure of $\pi$-determinism 
can be blamed on lack of communication among fibers. An example
suggested by Walters (see \cite{Walters}) also shows that there can be
{\em fully supported} $\nu$ on $Y$ which are 
not $\pi$-determinate. For such
examples there are potential functions 
$V \in \cc(Y)$ such that $V \circ\pi$ has
two equilibrium states which project to the {\em same} ergodic measure 
on $Y$.

In this example, $X=Y=\Sigma_2=$ full $2$-shift, and $\pi (x)_0=x_0+x_1 
\mod 2$ is a simple cellular automaton $2$-block map. If we replace
$X$ by its 
$2$-block recoding, so that $\pi$ becomes a $1$-block map, we obtain
the following diagram:

\bigskip
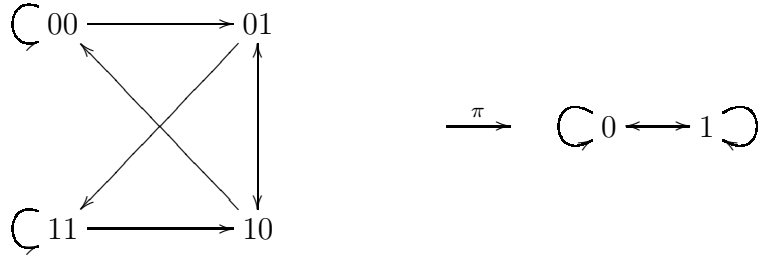
\begin{center}
\xymatrix{
& & 00 \ar@(ul,dl) \ar@{->}[rr] & &     01  \ar@{<->} [dd] \ar@{->}[ddll]\\
& & & & & & \ar [r]^\pi  & & 0 \ar@(ul,dl) \ar@{<->}
[r] & 1 \ar@(ur,dr)\\
& & 11 \ar@(ul,dl)  \ar@{->}[rr] & & 10 \ar@{->}[uull]
}
\end{center}

This is a finite-to-one map and hence is Markovian---for example, the
Bernoulli $1/2, 1/2$ measure on $\Sigma_2$ is mapped to itself. The
constant function $0$ is a compensation function. Thus every Markov
measure on $Y$ is $\pi$-determinate: the equilibrium state $\mu_V$ of a
locally constant $V$ on $Y$ lifts to the equilibrium state of $V \circ 
\pi$, which is the unique relatively maximal measure over $\mu_V$ (in
fact it's the only measure in $\pi^{-1}\{\mu_V\}$).

For every ergodic $\nu$ on $Y$,
{\em all} of $\pi^{-1}\{\nu\}$ consists of relatively maximal measures
over $\nu$, all of them having the same entropy as $\nu$.

If $p \neq 1/2$, the two measures on the SFT $X$ that
correspond to the Bernoulli measures $\cb (p,1-p)$ and $\cb (1-p,p)$
both map to the same measure $\nu_p$ on $Y$. Thus $\nu_p$, which is
fully supported on $Y$, is {\em not} $\pi$-determinate.
(An entropy-decreasing example is easily produced by forming the
Cartesian product of $X$ with another SFT.)

Moreover, $\nu_p$ is the unique equilibrium state of some continuous
function $V_p$ on $Y$ \cite{Phelps}. Then the set of relatively maximal measures over
$\nu_p$, which is the entire set $\pi^{-1}\{\nu_p\}$, consists of 
 the equilibrium
states of $V_p \circ\pi + G \circ \pi = V_p \circ \pi$ \cite{Walters}, so this
potential function $V_p \circ \pi$ has many equilibrium states.
\end{ex}


\begin{ex}\label{homclumps}
{\em Homogeneous clumps}.
In the following example there is no singleton clump, but the clumps
are homogeneous with respect to $\pi$ so there is a locally constant
compensation function (see \cite{Boyle-T, Shin, Shin2}), 
and hence every Markov measure on $Y$ is
$\pi$-de\-term\-in\-ate and its unique relatively maximal lift is Markov.

\bigskip
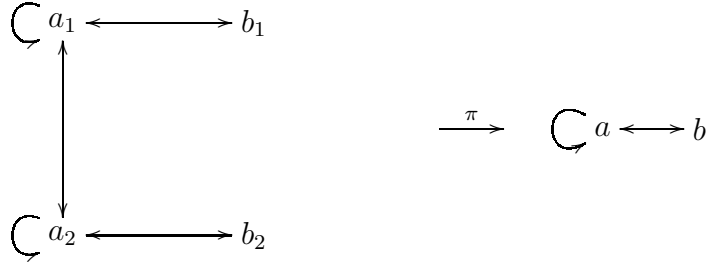
\begin{center}
\xymatrix{
& & a_1 \ar@(ul,dl) \ar@{<->}[dd] \ar@{<->}[rr] & &     b_1  \\
& & & & & & \ar [r]^\pi & & a \ar@(ul,dl) \ar@{<->}
[r] & b \\
& & a_2 \ar@(ul,dl) \ar@{<->}[rr] & & b_2 
}
\end{center}

\bigskip
\noindent
In this case the return time to $[a]$ is bounded, so $X_a$ is a
finite-state SFT rather than the countable-state chain of the general
case. There are six states, $a_1a_1$, 
$a_1b_1a_1$, $a_1a_2$, $a_2a_2$, $a_2b_2a_2$, and $a_2a_1$, 
according to the time 0
entries of $x \in X_a$ and $\sigma_ax$. Fix this order of the states
for indexing purposes. It can be shown by direct calculation that for
this example a
stochastic matrix $P$ determines a Markov measure on $X_a$ that is
relatively maximal over its image if and only if it is of the form
\begin{equation}\label{P}
\left(
\begin{matrix}
  x & 1-2x & x & 0 & 0 & 0 \\
  y & 1-2y & y & 0 & 0 & 0 \\
  0 & 0 & 0 & x & 1-2x & x \\
  0 & 0 & 0 & x & 1-2x & x \\
  0 & 0 & 0 & y & 1-2y & y \\
  x & 1-2x & x & 0 & 0 & 0
\end{matrix} 
\right) \quad .
\end{equation}
(In this case the image measure is also Markov.)\\
Here $0<x,y<1/2$ and the probability vector fixed by $P$ is 
\begin{equation}
p=\frac{1}{4y+2(1-2x)}(y,1-2x,y,y,1-2x,y) .
\end{equation}
Further, given a (1-step) Markov measure $\nu$ on $Y$, put $K=\nu
[aa]/\nu [aba]$. Then a stochastic matrix of the form (\ref{P})
with fixed vector $p$ satisfies $p_1+p_3+p_4+p_6=\nu [aa]$ and
$p_2+p_5=\nu [aba]$ (so that the Markov measure $\mu$ that it determines
projects to $\nu$) if and only if $x=y=K/(2K+2)$ (and then $\mu$ is
relatively maximal over $\nu$). 
\end{ex}

\begin{ex}\label{higher}
{\em Singleton clump after recoding}.
Make the preceding example a little bit more complicated by adding a
loop at $b_1$, so that now the return time to $[a]$ is unbounded. It
can be verified that now there is still a continuous saturated
compensation function, but there is no locally constant compensation
function, so the code is not Markovian. 
However, if we look at higher block presentations of $X$ and
$Y$, we can find singleton clumps, for example $abba$. Therefore again 
every Markov measure on $Y$ is $\pi$-determinate.

\bigskip
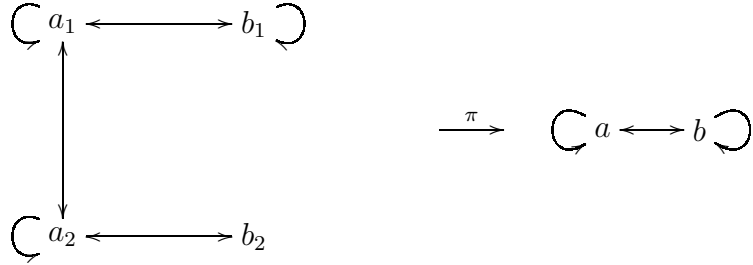
\begin{center}
\xymatrix{
& & a_1 \ar@(ul,dl) \ar@{<->}[dd] \ar@{<->}[rr] & &     b_1 \ar@(ur,dr) \\
& & & & & & \ar [r]^\pi & & a \ar@(ul,dl) \ar@{<->}
[r] & b \ar@(ur,dr)\\
& & a_2 \ar@(ul,dl) \ar@{<->}[rr] & & b_2 
}
\end{center}
\end{ex}

\begin{ex}\label{unknown}
{\em No singleton clumps}.
Complicating Example \ref{higher} a bit more, we can produce a situation in
which there are no singleton clumps, not even for any higher block
presentation.

\bigskip
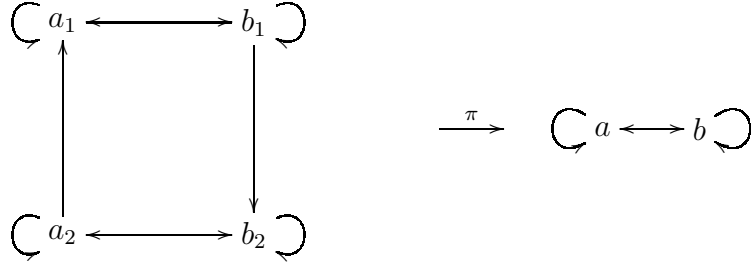
\begin{center}
\xymatrix{
& & a_1 \ar@(ul,dl) \ar@{<->}[rr] & &     b_1 \ar@(ur,dr) \ar 
[dd]\\
& & & & & & \ar [r]^\pi  & & a \ar@(ul,dl) \ar@{<->}
[r] & b \ar@(ur,dr)\\
& & a_2 \ar@(ul,dl) \ar [uu] \ar@{<->}[rr] & & b_2 \ar@(ur,dr) 
}
\end{center}

\bigskip
\noindent
For this example it can be shown that there is a continuous saturated 
compensation function $G \circ \pi$, but we do not know 
exactly which measures are $\pi$-determinate. Although the 
example appears simple, the question of how many fibers allow how much 
switching is complex. 
\end{ex}

\begin{ack*} An exposition of some of these results was included in
  \cite{Petersen00}. We thank Mike Boyle, Xavier M\' ela,  Jean-Paul
  Thouvenot, and Peter Walters for their
  very helpful
  questions, suggestions, and discussions and the referee for comments
  that improved the writing.
\end{ack*}
\bibliographystyle{amsplain}
\providecommand{\bysame}{\leavevmode\hbox to3em{\hrulefill}\thinspace}

\end{document}